\documentclass[11pt]{article}
\usepackage{amsmath,latexsym}
\usepackage{enumerate}
\usepackage{epsfig}

\newcommand\RR{\rm {I \! R}}

\newcommand\NN{\rm {I \! N}}

\newcommand\TT{{\cal T}}

\DeclareMathOperator\sP{P}   
\newcommand{\rP}{\mathrm{P}} 
\DeclareMathOperator\sE{E}   
\newcommand{\rE}{\mathrm{E}} 
\DeclareMathOperator\sI{I}   
 
\newcommand{\eps}{\varepsilon}

\newcommand\address{{\bf Author's address}:
Mathematics Department, University of North Texas,
P.O. Box 311430, Denton, TX 76203-1430; e-mail: allaart@unt.edu}

\newtheorem{theorem}{Theorem}[section]

\newtheorem{corollary}[theorem]{Corollary}
\newtheorem{example}[theorem]{Example}
\newtheorem{lemma}[theorem]{Lemma}
\newtheorem{remark}[theorem]{Remark}

\title{Prophet inequalities for i.i.d. random variables \\with random arrival times.}
\author{Pieter C. Allaart\\University of North Texas \footnote{\address}}
\date{\today}

\begin{document}

\maketitle

\begin{abstract}
Suppose $X_1,X_2,\dots$ are i.i.d. nonnegative random variables with finite expectation, and for each $k$, $X_k$ is observed at the $k$-th arrival time $S_k$ of a Poisson process with unit rate which is independent of the sequence $\{X_k\}$. For $t>0$, comparisons are made between the expected maximum $M(t):=\rE[\max_{k\geq 1} X_k \sI(S_k\leq t)]$ and the optimal stopping value $V(t):=\sup_{\tau\in\TT}\sE[X_\tau \sI(S_\tau\leq t)]$, where $\TT$ is the set of all $\NN$-valued random variables $\tau$ such that $\{\tau=i\}$ is measurable with respect to the $\sigma$-algebra generated by $(X_1,S_1),\dots,(X_i,S_i)$. For instance, it is shown that $M(t)/V(t)\leq 1+\alpha_0$, where $\alpha_0\doteq 0.34149$ satisfies $\int_0^1(y-y\ln y+\alpha_0)^{-1}\,dy=1$; and this bound is asymptotically sharp as $t\to\infty$. Another result is that $M(t)/V(t)<2-(1-e^{-t})/t$, and this bound is asymptotically sharp as $t\downarrow 0$. Upper bounds for the difference $M(t)-V(t)$ are also given, under the additional assumption that the $X_k$ are bounded.
\end{abstract}



\section{Introduction}

Suppose an item is for sale, and independent, identically distributed price offers $X_1,X_2,\dots$ arrive according to a Poisson process with rate $\lambda$. The item must be sold by a fixed time $t>0$, or it becomes worthless. If $S_i$ denotes the arrival time of the $i$-th offer, the optimal expected return is given by
\begin{equation*}
V(t):=\sup_{\tau\in\TT}\sE[X_\tau \sI(S_\tau\leq t)],
\end{equation*}
where $\TT$ is the set of all $\NN$-valued random variables (stopping rules) $\tau$ such that $\{\tau=i\}$ is measurable with respect to the $\sigma$-algebra generated by $(X_1,S_1),\dots,(X_i,S_i)$. 
Assuming the price offers are independent of the arrival process, Karlin (1962) showed that $V(t)$ is the unique solution of the initial value problem
\begin{equation}
V'=\lambda\sE(X_1-V)^+, \qquad V(0)=0,
\label{eq:V-diff-eq}
\end{equation}
and the optimal policy is to accept the first offer whose value $X$ exceeds $V(\tau)$, where $\tau$ is the amount of time remaining when the offer arrives. Sakaguchi (1976) gives explicit solutions of \eqref{eq:V-diff-eq} for several common distributions of $X_1$.

The purpose of this paper is to compare $V(t)$ with the expected maximum
\begin{equation*}
M(t):=\sE(\max\{X_1,\dots,X_{N(t)}\}),
\end{equation*}
where $N(t)$ denotes the number of arrivals up to time $t$. (The maximum of an empty set is taken to be zero.) In particular, reasonably tight upper bounds are given for the ratio $M(t)/V(t)$ and the difference $M(t)-V(t)$, the latter under the additional assumption that the price offers are uniformly bounded. 

Note that $M(t)$ may be interpreted as the optimal expected return to a player endowed with complete foresight, who is able to predict both the values and the arrival times of future price offers. Such a player would of course accept the largest offer made available before the deadline. Hence, comparisons of $V(t)$ and $M(t)$ can be interpreted as bounds on the value of inside information in investing. For instance, consider two investors holding identical assets. One investor has no inside information, and therefore does not know when offers will be made or what the sizes of the offers will be. The other, better connected, investor happens to know that a certain firm is considering to make a bid within the next few days, and perhaps even knows the likely size of the bid. The inequalities developed in this paper can be used to bound the advantage of the second investor over the first investor.

Analogous comparisons in the discrete-time setting are known in the literature as {\em prophet inequalities}. Specifically, let $Z_1,\dots,Z_n$ be independent nonnegative random variables with finite expectations, and define the quantities
\begin{equation*}
V_n:=V(Z_1,\dots,Z_n):=\sup\{\sE Z_\tau:\tau\ \mbox{is a stopping rule for}\ Z_1,\dots,Z_n\}, 
\end{equation*}
and
\begin{equation*}
M_n:=\sE(\max\{Z_1,\dots,Z_n\}). 
\end{equation*}
Krengel and Sucheston (1978) proved that $M_n\leq 2V_n$, with strict inequality if $Z_1,\dots,Z_n$ are not all identically equal to zero. Hill and Kertz (1981) showed that if, instead of being nonnegative, the $Z_i$ take values in an interval $[a,b]$, then $M_n-V_n\leq (b-a)/4$. While the constants $2$ and $1/4$ are best-possible in general, sharper bounds hold if the $Z_i$ are known to be identically distributed. Hill and Kertz (1982) constructed the best-possible constants $a_n$ and $b_n$ ($n=2,3,\dots$) such that {\em (a)} if $Z_1,\dots,Z_n$ are i.i.d. nonnegative random variables, then $M_n\leq a_n V_n$; and {\em (b)} if $Z_1,\dots,Z_n$ are i.i.d. taking values in $[a,b]$, then $M_n-V_n\leq b_n(b-a)$. The constants $a_n$ and $b_n$ play a prominent role in the present article. Their rather technical definitions are given in Section \ref{notation} below.

There are a few important differences between the classical discrete-time setting and the continuous-time setting considered here. First, since there is no upper bound on the possible number of observations in the interval $[0,t]$, the usual technique of backward induction can not be applied here. This makes it impossible to reduce the problem of finding sharp inequalities to a finite-dimensional optimization problem, as was done by Hill and Kertz (1982). Second, even though $X_1,X_2,\dots$ are i.i.d., the random variables $X_k \sI(S_k\leq t), k\in\NN$ are in fact positively dependent, with expectations decreasing to zero.

Given these difficulties, it seems unlikely that completely sharp but explicit prophet inequalities can be found for the continuous-time setting. In this paper, four inequalities are proved: two for the ratio $M(t)/V(t)$, and two for the difference $M(t)-V(t)$, the latter under the additional assumption that the $X_k$ are bounded. The upper bounds given in Section~\ref{long-range} are constant in $t$, and are asymptotically sharp as $t\to\infty$. These bounds make use of the Hill-Kertz constants $\{a_n\}$ and $\{b_n\}$. Obviously, such bounds are far from sharp when $t$ is small, since the prophet's advantage disappears as the expected number of observations approaches zero. Therefore, a second pair of bounds, which are asymptotically sharp as $t\downarrow 0$, is given in Section~\ref{short-range}. The proofs of these bounds rely on the use of threshold stopping rules.

Section \ref{examples} provides some estimates on how close to being sharp the bounds of Sections~\ref{long-range} and \ref{short-range} are. Finally, Section~\ref{renewal} discusses more general renewal processes.

Throughout the paper it is assumed that $\lambda=1$. However, all of the results can be restated easily (with only trivial modifications) for arbitrary $\lambda>0$ and, more generally, for any non-homogeneous Poisson process.

\section{Notation and the Hill-Kertz constants} \label{notation}

Throughout the paper, $X,X_1,X_2,\dots$ denote nonnegative i.i.d. random variables with finite expectation, where it is assumed that $\rP(X=0)<1$. For real numbers $x$ and $y$, $x\vee y$ denotes the maximum of $x$ and $y$.  To emphasize the dependence of $M(t)$ and $V(t)$ on the distribution of $X$, they will sometimes be written as $M(t;X)$ and $V(t;X)$, respectively. For $s>0$, set $X_s^*:=\max\{X_1,\dots,X_{N(s)}\}$.

We next introduce the constants $\{a_n\}$ and $\{b_n\}$ mentioned in the introduction. The definitions below are taken from Hill and Kerz (1982). For $n>1$ and $w,x\in[0,\infty)$, let $\phi_n(w,x)=(n/(n-1))w^{(n-1)/n}+x/(n-1)$. For $\alpha\in[0,\infty)$, define $\eta_{0,n}(\alpha)=\phi_n(0,\alpha)$, and inductively, $\eta_{j,n}(\alpha)=\phi_n(\eta_{j-1,n}(\alpha),\alpha)$ for $j\geq 1$. In their Propositions 3.4 and 3.8, Hill and Kertz show that

{\em (i)} there is a unique $\alpha_n\in(0,1)$ such that $\eta_{n-1,n}(\alpha_n)=1$; and

{\em (ii)} there is a unique $\beta_n\in(0,1)$ such that $(n-1)[\eta_{n,n}(\beta_n)-\eta_{n-1,n}(\beta_n)]=1$.

\noindent Kertz (1986; Lemma 6.2(b)) proves further that $\lim_{n\to\infty}\alpha_n=\alpha_0$, where $\alpha_0\doteq 0.34149$ is the unique value of $\alpha$ such that $\int_0^1(y-y\ln y+\alpha)^{-1}\,dy=1$. Table \ref{tab:alphas} below gives sample values of $\alpha_n$ and $\beta_n$.

\begin{table}[h]
\[
\begin{array}{rcc}
n & \alpha_n & \beta_n \\ \hline
2 & .17157 & .06250\\
3 & .22138 & .07761\\
4 & .24811 & .08539\\
5 & .26496 & .09020\\
6 & .27659 & .09348\\
7 & .28513 & .09586
\end{array}
\hskip 1in
\begin{array}{rcc}
n & \alpha_n & \beta_n \\ \hline
8 & .29166 & .09768\\
9 & .29683 & .09911\\
10 & .30101 & .10027\\
100 & .33716 & .11010\\
10^4 & .34144 & .11125\\
10^6 & .34149 & .11126
\end{array}
\]
\caption{Selected values of $\alpha_n$ and $\beta_n$.}
\label{tab:alphas}
\end{table}

Now let $a_n:=1+\alpha_n$, and $b_n:=\beta_n$. 

\begin{theorem} \label{thm:HK}
(a) {\em [Hill and Kertz, 1982; Theorem A]} If $Z_1,\dots,Z_n$ are i.i.d. nonnegative random variables, then $\rE(Z_1\vee\dots\vee Z_n)\leq a_n V(Z_1,\dots,Z_n)$. This bound is sharp, and holds with strict inequality if $Z_1$ is not identically equal to zero.

(b) {\em [Hill and Kertz, 1982; Theorem B]} If $Z_1,\dots,Z_n$ are i.i.d. random variables taking values in $[a,b]$, then  $\rE(Z_1\vee\dots\vee Z_n)-V(Z_1,\dots,Z_n)\leq b_n(b-a)$. This bound is attained.
\end{theorem}

\section{Long-range prophet inequalities} \label{long-range}

This section develops upper bounds for the ratio and difference of $M(t)$ and $V(t)$ which are fairly sharp when $t$ is large. The following simple inequality will be helpful.

\begin{lemma} \label{simple-bound}
For all $s>0$ and all $c\geq 0$,
\begin{equation*}
\rE(X_s^*-c)^+\leq s\sE(X-c)^+,
\end{equation*}
with strict inequality if $\rP(X>c)>0$.
\end{lemma}

\noindent {\em Proof.} This follows easily by conditioning on $N(s)$, and using the fact that $\sE(X_1\vee\dots\vee X_n-c)^+\leq n\sE(X-c)^+$ for every $n$, with strict inequality if $n\geq 2$ and $\rP(X>c)>0$. $\Box$

\begin{theorem} \label{thm:uniform-ratio-bound}
Let $\alpha_0$ be as in Section \ref{notation}. For all $t>0$,
\begin{equation*}
M(t)\leq (1+\alpha_0)V(t).
\end{equation*}
\end{theorem}

\noindent {\em Proof.} Fix $n\geq 2$, and let $\delta:=t/n$. Consider a ``partial prophet" who has limited foresight in the sense that he can see, at the beginning of each time interval $I_i:=((i-1)\delta,i\delta]$, $i=1,\dots,n$, all of the observations (if any) arriving in that interval. For $i=1,\dots,n$, let $Z_i$ denote the largest of the observations arriving in the interval $I_i$ (or $Z_i=0$ if no observations arrive in that interval). A routine exercise shows that $Z_1,\dots,Z_n$ are i.i.d.

Let $v_j:=V(Z_{j+1},\dots,Z_n)$, $j=0,1,\dots,n-1$, and put $v_n:=0$. Since the partial prophet sequentially observes  $Z_1,\dots,Z_n$, his optimal expected return is $v_0$, and by backward induction (see p.~50 of Chow et al., 1971), his optimal rule is to stop in the first time interval $I_i$ for which $Z_i\geq v_i$, and to accept the largest observation, $Z_i$, in that interval.

Now consider the following stopping rule for the mortal:
\begin{quote}
Accept the first observation $X_j$ such that, if $X_j$ arrives in the time interval $I_i$, then $X_j\geq v_i$.
\end{quote}
Let $V_\delta(t)$ denote the expected return from this stopping rule. Define $Z_i'=Z_i\sI(Z_i\geq v_i)$, and let 
$X_{i,1}',X_{i,2}',\dots$ denote the successive values of those $X_j$ arriving after time $(i-1)\delta$ for which $X_j\geq v_i$. Let $N_i'$ be the number of such observations (with $X_j\geq v_i$) that arrive in the interval $I_i$. Observe that 
$N_i'$ is Poisson with mean $\delta\sP(X\geq v_i)$, and $Z_i'=\max\{X_{i,1}',\dots,X_{i,N_i'}'\}$. 
Thus, Lemma \ref{simple-bound} applied to $X_{i,1}',X_{i,2}',\dots$ yields
\begin{equation}
\rE Z_i'\leq\delta\sP(X\geq v_i)\sE X_{i,1}'.
\label{eq:ZX-relationship}
\end{equation}
Note that
\begin{equation}
v_0=\sum_{i=1}^n\sP(Z_1<v_1,\dots,Z_{i-1}<v_{i-1})\sE[Z_i\sI(Z_i\geq v_i)],
\label{eq:partial-prophet-value}
\end{equation}
and
\begin{equation}
V_\delta(t)=\sum_{i=1}^n\sP(Z_1<v_1,\dots,Z_{i-1}<v_{i-1})\sE[X_{i,1}'\sI(Z_i\geq v_i)].
\label{eq:gambler-lower-bound}
\end{equation}
Using \eqref{eq:ZX-relationship}, we obtain that
\begin{align*}
\sE[X_{i,1}'\sI(Z_i\geq v_i)]&=\sE X_{i,1}'\sP(N_i'\geq 1)
\geq \sE Z_i'\cdot \frac{\rP(N_i'\geq 1)}{\delta\sP(X\geq v_i)}\\
&=\sE Z_i'\cdot \frac{1-e^{-\delta\sP(X\geq v_i)}}{\delta\rP(X\geq v_i)}
\geq \sE Z_i'\cdot(1-e^{-\delta})/\delta\\
&=\sE[Z_i\sI(Z_i\geq v_i)](1-e^{-\delta})/\delta,
\end{align*}
where the second inequality follows since $(1-e^{-\delta p})/p$ is decreasing in $p$. Substituting this result into \eqref{eq:gambler-lower-bound} and comparing with \eqref{eq:partial-prophet-value} yields the conclusion
\begin{equation}
v_0\leq \frac{\delta}{1-e^{-\delta}} V_\delta(t)\leq\frac{t/n}{1-e^{-t/n}}V(t).
\label{eq:partial-prophet-inequality}
\end{equation}
By Theorem \ref{thm:HK}(a),
\begin{equation*}
M(t)=\sE(X_1\vee\dots\vee X_{N(t)})=\sE(Z_1\vee\dots\vee Z_n)<a_n v_0.
\end{equation*}
Hence,
\begin{equation*}
M(t)<a_n \frac{t/n}{1-e^{-t/n}}V(t).
\end{equation*}
Since $n$ was arbitrary, the theorem follows upon letting $n\to\infty$. $\Box$

\begin{theorem} \label{thm:uniform-difference-bound}
Assume that $X_1,X_2,\dots$ are $[0,1]$-valued. Then for all $t>0$,
\begin{equation*}
M(t)-V(t)\leq \limsup_{n\to\infty} b_n.
\end{equation*}
\end{theorem}

\noindent {\em Proof.} We use the notation from the proof of Theorem \ref{thm:uniform-ratio-bound}. Since $v_0\leq \sP(N(t)\geq 1)=1-e^{-t}$, Theorem \ref{thm:HK}(b) and \eqref{eq:partial-prophet-inequality} imply that
\begin{equation}
M(t)-V(t)\leq b_n+(1-e^{-t})\left[1-\frac{n(1-e^{-t/n})}{t}\right], \qquad\mbox{for all $n\geq 2.$}
\label{eq:precise-bound}
\end{equation}
Letting $n\to\infty$ completes the proof. $\Box$

\begin{remark}
{\rm
The value of $\limsup b_n$ does not seem to be known at present, though Table~\ref{tab:alphas} ($b_n=\beta_n$) suggests that $\limsup b_n\approx .1113$. Fortunately, equation \eqref{eq:precise-bound} gives rigorous upper bounds by taking $n=10^6$. For instance, if $t=1000$, \eqref{eq:precise-bound} yields that $M(t)-V(t)\leq .11176$.
}
\end{remark}

\begin{remark}
{\rm
Theorems \ref{thm:uniform-ratio-bound} and \ref{thm:uniform-difference-bound} hold in fact for any non-homogeneous Poisson process with rate function $\lambda(t)$, $t>0$, provided $\lambda(t)$ is bounded on bounded intervals. This follows since in the proof of Theorem \ref{thm:uniform-ratio-bound}, one may replace the partition $\{I_i\}$ with a partition in which the number of arrivals in each interval is Poisson with the {\em same} parameter $\mu=(1/n)\int_0^t\lambda(s)\,ds$. The rest of the proof then goes through with $\mu$ in place of $\delta$.
}
\end{remark}

\begin{remark}
{\rm
It is possible to quantify how sharp the bounds of Theorems \ref{thm:uniform-ratio-bound} and \ref{thm:uniform-difference-bound} are: Proposition 4.4 of Hill and Kertz (1982) gives the $\eps$-extremal distributions for Theorem \ref{thm:HK}(a). They satisfy $\rP(Z_1=0)=(\eta_{0,n}(\alpha_n))^{1/n}=(\alpha_n/(n-1))^{1/n}$. By the construction of the random variables $\{Z_i\}$ in the proof of Theorem \ref{thm:uniform-ratio-bound}, $Z_1$ must have an atom at zero of size at least $\rP(N(\delta)=0)=e^{-\delta}=e^{-t/n}$. Vice versa, every distribution on $[0,\infty)$ satisfying this condition can arise from a suitable choice of the distribution of $X$. It follows that, if
\begin{equation*}
t\geq \log((n-1)/\alpha_n),
\end{equation*}
then for every $\eps>0$ there exists a random variable $X$ and corresponding i.i.d. random variables $Z_1,\dots,Z_n$ such that
\begin{equation*}
M(t;X)=\sE(Z_1\vee\dots\vee Z_n)>(a_n-\eps)V(Z_1,\dots,Z_n)\geq (a_n-\eps)V(t;X).
\end{equation*}
For example, if $t\geq \log(99/\alpha_{100})\doteq 5.683$, then $\sup_X M(t;X)/V(t;X)\geq a_{100}\doteq 1.337$.

Similarly, the extremal distribution for Theorem \ref{thm:HK}(b) (see Proposition 5.3 of Hill and Kertz, 1982) has $\rP(Z_1=0)=(\eta_{0,n}(\beta_n))^{1/n}$, so as long as $t\geq\log((n-1)/\beta_n)$, there exists a $[0,1]$-valued random variable $X$ such that $M(t;X)-V(t;X)\geq b_n$. For example, if $t\geq \log(99/\beta_{100})\doteq 6.802$, then $M(t;X)-V(t;X)\geq b_{100}\doteq 0.110$ for a suitable $X$ in $[0,1]$.
}
\end{remark}

\begin{remark}
{\rm
In contrast with the classical discrete-time setting, there is no obvious generalization of Theorem~\ref{thm:uniform-difference-bound} to random variables taking values in an arbitrary interval $[a,b]$. Indeed, if $\hat{X}_k$ is in $[a,b]$, the random variables $\hat{X}_k \sI(S_k\leq t)$ are not in $[a,b]$ but in $\{0\}\cup[a,b]$. Hence the difference $M(t)-V(t)$ is not invariant under a shift of the distribution of $X$. A simple example illustrates this: Let $X$ take the values $0$ and $1$ with probabilities $1/2$ each; clearly, $M(t;X)=V(t;X)$. Now set $\hat{X}:=X+1$, and consider the problem of stopping the sequence $\hat{X}_k \sI(S_k\leq t)$. If $t$ is  sufficiently small, it is optimal for the mortal to accept the first available observation regardless of its value, whereas the prophet can wait for a larger value (which will arrive with positive probability). Hence $V(t;\hat{X})=\sE[\hat{X}_1 \sI(S_1\leq t)]=(1-e^{-t})\sE\hat{X}<M(t;\hat{X})$.

On the other hand, if $X$ takes values in $[a,b]$ with $0\leq a<b$, then the proofs of Theorems \ref{thm:uniform-ratio-bound} and \ref{thm:uniform-difference-bound} are easily modified to show that $M(t)-V(t)\leq b\cdot\limsup_{n\to\infty}b_n$. This bound is not, however, very sharp in general.
}
\end{remark}

\section{Short-range prophet inequalities} \label{short-range}

The inequalities obtained in the previous section can be improved upon when $t$ (and with it the expected number of observations) is small. In this section, we consider pure threshold rules of the form $\tau(c)=\inf\{n:X_n\geq c\}$.  For $c\geq 0$, let $W_c(t)$ denote the expected return from the rule $\tau(c)$. That is,
\begin{equation*}
W_c(t):=W_c(t;X):=\sE\left[X_{\tau(c)}\sI(\tau(c)\leq N(t))\right].
\end{equation*}
It is straightforward to verify that the value of $W_c(t)$ is given by
\begin{equation}
W_c(t)=\left[1-e^{-t\sP(X\geq c)}\right]\sE(X|X\geq c).
\label{eq:threshold-value}
\end{equation}

The next lemma is the key to the results in this section.

\begin{lemma} \label{lem:calculus}
For $\gamma>0$, the function
\begin{equation*}
f(x):=(1-e^{-x})(1+\gamma/x), \qquad x>0
\end{equation*}
does not have a local minimum on $(0,\infty)$.
\end{lemma}

\noindent {\em Proof.} Twice differentiating $f$ yields that
\begin{equation}
f'(x)+f''(x)=(\gamma/x^3)\{2-x-(2+x)e^{-x}\}=:(\gamma/x^3)g(x).
\label{eq:derivsum}
\end{equation}
Observe that $g(0)=g'(0)=0$, and $g''(x)<0$ for all $x>0$. Thus, by \eqref{eq:derivsum}, $f'(x)+f''(x)<0$ for all $x>0$. But then there cannot exist a point $x_0>0$ such that $f'(x_0)=0$ and $f''(x_0)\geq 0$. Since $f$ is smooth, the lemma follows. $\Box$

\begin{theorem} \label{thm:threshold-bound}
For all $t>0$,
\begin{equation}
\frac{M(t)}{\sup_c W_c(t)}<2-\frac{1-e^{-t}}{t},
\label{eq:threshold-bound}
\end{equation}
and this bound is sharp.
\end{theorem}

\begin{corollary} \label{cor:short-range-ratio-bound}
For all $t>0$,
\begin{equation*}
\frac{M(t)}{V(t)}<2-\frac{1-e^{-t}}{t}.
\end{equation*}
\end{corollary}

\noindent {\em Proof of Theorem \ref{thm:threshold-bound}.} For any $c\geq 0$,
\begin{align*}
\rE[X_t^*\sI(X_t^*<c)]&\leq c\sP(X_t^*<c, N(t)\geq 1)\\
&=c\{\sP(X_t^*<c)-P(N(t)=0)\}\\
&=c\sP(X_t^*<c)-c e^{-t},
\end{align*}
and
\begin{equation*}
\rE[X_t^*\sI(X_t^*\geq c)]=c\sP(X_t^*\geq c)+\sE(X_t^*-c)^+.
\end{equation*}
Adding these expressions and applying Lemma \ref{simple-bound} gives
\begin{equation}
\rE X_t^*\leq c(1-e^{-t})+\sE(X_t^*-c)^+\leq c(1-e^{-t})+t\sE(X-c)^+.
\label{eq:simpler-max-upper-bound}
\end{equation}
On the other hand, \eqref{eq:threshold-value} can be written as
\begin{equation}
W_c(t)=\left[1-e^{-t\sP(X\geq c)}\right]\left\{c+\frac{\rE(X-c)^+}{\rP(X\geq c)}\right\}.
\label{eq:threshold-value2}
\end{equation}
Now for any $\alpha>0$, there exists a unique number $c:=c_\alpha$ such that
\begin{equation}
\rE(X-c)^+=c\alpha.
\label{eq:c-equation}
\end{equation}
For this value of $c$, \eqref{eq:simpler-max-upper-bound} and \eqref{eq:threshold-value2} combine to give
\begin{equation*}
\frac{M(t)}{W_c(t)}\leq \frac{1-e^{-t}+t\alpha}{(1-e^{-tp})(1+\alpha/p)},
\end{equation*}
where $p:=\sP(X\geq c)$. Let $d(p):=(1-e^{-tp})(1+\alpha/p)$. Using Lemma \ref{lem:calculus} with $x=tp$ and $\gamma=t\alpha$ we find that $d(p)$ is smallest either at $p=0^+$ or at $p=1$. Thus, noting that $d(0^+)=t\alpha$,
\begin{equation}
\frac{M(t)}{W_c(t)}\leq \max\left\{\frac{1-e^{-t}+t\alpha}{t\alpha},\frac{1-e^{-t}+t\alpha}{(1-e^{-t})(1+\alpha)}\right\}.
\label{eq:uniform-ratio-bound}
\end{equation}
The first term in the maximum is decreasing, and the second increasing in $\alpha$. Hence, the right hand side of \eqref{eq:uniform-ratio-bound} is minimized when $t\alpha=(1-e^{-t})(1+\alpha)$; that is, when
\begin{equation*}
\alpha=\alpha^*:=\frac{1-e^{-t}}{t+e^{-t}-1}.
\end{equation*}
(Note that $\alpha^*>0$.) For $\alpha=\alpha^*$, the maximum in \eqref{eq:uniform-ratio-bound} reduces to the right hand side of \eqref{eq:threshold-bound}. Finally, the inequality is strict since \eqref{eq:c-equation} and $\alpha^*>0$ imply that $\rP(X>c_{\alpha^*})>0$, giving strict inequality in Lemma \ref{simple-bound} (and hence, in \eqref{eq:simpler-max-upper-bound}).

To see that the bound is sharp, let $0<p<1$, $\eps=\{1-e^{-pt}-p(1-e^{-t})\}/(1-p)(1-e^{-t})$, and let $X$ have the distribution $\rP(X=1)=p=1-\sP(X=\eps)$. There are only two essentially different threshold rules: $\tau(\eps)$ and $\tau(1)$. By the choice of $\eps$,
\begin{equation*}
W_\eps(t)=(1-e^{-t})\{p+\eps(1-p)\}=1-e^{-tp}=W_1(t),
\end{equation*}
using \eqref{eq:threshold-value}. Furthermore, it is not difficult to compute that
\begin{equation*}
M(t)=1-(1-\eps)e^{-tp}-\eps e^{-t}.
\end{equation*}
Thus,
\begin{align*}
\frac{M(t)}{\sup_c W_c(t)}&=\frac{M(t)}{W_1(t)}
=1+\frac{\eps(e^{-tp}-e^{-t})}{1-e^{-tp}}\\
&=1+\frac{e^{-tp}-e^{-t}}{(1-p)(1-e^{-t})}\left(1-\frac{p(1-e^{-t})}{1-e^{-tp}}\right)\\
&\to 2-\frac{1-e^{-t}}{t} \qquad\mbox{as $p\downarrow 0$}.
\end{align*}
This shows that the inequality \eqref{eq:threshold-bound} is sharp. $\Box$

\bigskip
The next result is a difference inequality for the case when $X$ is bounded. Some additional notation is needed. For constants $a<b$, let $\mathcal{X}_{[a,b]}$ denote the collection of all $[a,b]$-valued random variables. For any threshold $c$, define $D_c(t;X):=M(t;X)-W_c(t;X)$. For $t>0$, let $h_t:[0,1]\to\RR$ be the function
\begin{equation*}
h_t(x)=1-e^{-tx}-(1-e^{-t})x.
\end{equation*}
Define $\gamma(t):=t/(1-e^{-t})$, and $\beta(t):=1-\{1+\log\gamma(t)\}/\gamma(t)$. Routine calculus shows that $\beta(t)=\max_{0\leq x\leq 1}h_t(x)$, and the maximum is attained at $x=(\log \gamma(t))/t$.

The next minimax-type theorem presents a universal value $c^*$ which minimizes the largest possible difference $D_c(t;X)$ for $X\in\mathcal{X}_{[a,b]}$.

\begin{theorem} \label{thm:minimax-bound}
Let $0\leq a<b$. For all $t>0$,
\begin{equation*}
\inf_c \sup_{X\in\mathcal{X}_{[a,b]}} D_c(t;X)=[b-\max\{a,c^*\}]\beta(t),
\end{equation*}
where $c^*=b\beta(t)/\{\beta(t)+1-e^{-t}\}$. Moreover, the infimum is attained by the choice $c=c^*$.
\end{theorem}

The proof of Theorem \ref{thm:minimax-bound} uses the concept of balayage. Given $X\in \mathcal{X}_{[a,b]}$ and constants $a\leq c<d\leq b$, let $X_c^d$ denote a random variable such that $X_c^d=X$ if $X\not\in [c,d]$, $X_c^d=c$ with probability $(d-c)^{-1}\sE[(d-X)\sI(c\leq X\leq d)]$, and $X_c^d=d$ otherwise. It follows immediately that $\rE X_c^d=\sE X$, $\rE(X_c^d|X_c^d\geq c)=\sE(X|X\geq c)$, and $\rP(X_c^d\geq c)=\sP(X\geq c)$. Moreover, Lemma 2.2 of Hill and Kertz (1981) implies that if $Y$ is a random variable independent of both $X$ and $X_c^d$, then $\rE(X_c^d\vee Y)\geq\sE(X\vee Y)$.

\bigskip
\noindent {\em Proof of Theorem \ref{thm:minimax-bound}.} Assume first that $c^*>a$. Let $X$ be any $[a,b]$-valued random variable; choose $c\in[a,b]$, and define $\hat{X}=X_c^b$. By \eqref{eq:threshold-value}, $W_c(t;\hat{X})=W_c(t;X)$, and by the last remark in the previous paragraph, $M(t;\hat{X})\geq M(t;X)$. Therefore, $D_c(t;\hat{X})\geq D_c(t;X)$. So, by replacing $X$ with $\hat{X}$ if necessary, we may assume that $\rP(c<X<b)=0$.

Next, define $p:=\sP(X=b)$, and $r:=\sP(X\geq c)$.
For $n\geq 1$,
\begin{align*}
\rE(X_1\vee\dots\vee X_n)&\leq c \sP(X_1\vee\dots\vee X_n\leq c)+b \sP(X_1\vee\dots\vee X_n=b)\\
&=c(1-p)^n+b\{1-(1-p)^n\}=b-(b-c)(1-p)^n.
\end{align*}
Thus,
\begin{align}
M(t;X)&\leq \sum_{n=1}^\infty \left\{b-(b-c)(1-p)^n\right\}\sP(N(t)=n)\notag\\
&=(b-c)(1-e^{-tp})+c(1-e^{-t}). \label{eq:M-bound}
\end{align}
On the other hand, by \eqref{eq:threshold-value},
\begin{equation*}
W_c(t;X)=(1-e^{-tr})\left(c+\frac{(b-c)p}{r}\right).
\end{equation*}
For fixed $p$, this expression is minimized either when $r=p$ or $r=1$, in view of Lemma \ref{lem:calculus}. It follows that
\begin{equation*}
W_c(t;X)\geq \min\{(1-e^{-t})(c+(b-c)p),b(1-e^{-tp})\}.
\end{equation*}
Subtracting from \eqref{eq:M-bound} and rearranging terms, we obtain that
\begin{align*}
D_c(t;X)&\leq \max\{(b-c)h_t(p),c(e^{-tp}-e^{-t})\}\\
&\leq \max\left\{(b-c)\beta(t),c(1-e^{-t})\right\}.
\end{align*}
The two terms inside the maximum are equal when $c=c^*$, and so
\begin{equation*}
D_{c^*}(t;X)\leq (b-c^*)\beta(t).
\end{equation*}

Suppose next that $c^*\leq a$. Then the preceding argument (with $r=\sP(X\geq a)=1$) yields that
\begin{equation*}
D_{c^*}(t;X)=D_a(t;X)\leq (b-a)\beta(t).
\end{equation*}

Conversely, for any $c\geq 0$ and $\eps>0$ the distribution of $X$ can be chosen so that $D_c(t;X)\geq [b-\max\{a,c^*\}]\beta(t)-\eps$: 

(i) If $c\leq a$, take $X\in\{a,b\}$ with $\rP(X=b)=(\log\gamma(t))/t$. Then $D_c(t;X)=(b-a)\beta(t)$.

(ii) If $a<c<c^*$, take $X\in\{c,b\}$ with $\rP(X=b)=(\log \gamma(t))/t$. Then $D_c(t;X)=(b-c)\beta(t)>(b-c^*)\beta(t)$.

(iii) If $c\geq c^*$ and $c>a$, take $X\equiv c-\delta$, where $0<\delta<\min\{\eps,c-a\}$. Then $D_c(t;X)=(c-\delta)(1-e^{-t})\geq c^*(1-e^{-t})-\delta\geq (b-c^*)\beta(t)-\eps$. 

Thus, the choice $c^*$ is minimax, and the theorem follows. $\Box$

\begin{corollary} \label{threshold-difference-bound} If $X$ is $[a,b]$-valued with $0\leq a<b$, then for all $t>0$,
\begin{equation}
M(t;X)-V(t;X)\leq\min\left\{(b-a)\beta(t),\frac{b\beta(t)(1-e^{-t})}{\beta(t)+1-e^{-t}}\right\}.
\label{eq:d-bound}
\end{equation}
\end{corollary}

\begin{remark}
{\rm
For a non-homogeneous Poisson process with rate function $\lambda(x)$, the bounds corresponding to Theorems \ref{thm:threshold-bound} and \ref{thm:minimax-bound} are obtained by replacing $t$ with $\mu(t):=\int_0^t\lambda(x)\,dx$.
}
\end{remark}

\section{How sharp are the bounds?} \label{examples}

In this section, the ratio and difference of $M(t;X)$ and $V(t;X)$ are examined for random variables $X$ taking only finitely many values, say $a_1<a_2<\dots<a_n$, where $a_1\geq 0$. Put $a_0=0$. For $i=0,1,\dots,n$, define $r_i:=\sP(X\geq a_i)$, $\mu_i:=\sE(X-a_i)^+$, and $E_i:=\sE(X|X\geq a_i)$. Observe that $\mu_n=0$, and recursively, for $k=n,n-1,\dots,1$,
\begin{equation*}
\mu_{k-1}=\sE(X-a_k)^++r_k(a_k-a_{k-1})=\mu_k+r_k(a_k-a_{k-1}).
\end{equation*}
A moment's reflection reveals that there are critical times $0<t_1^*<t_2^*<\dots<t_{n-1}^*<\infty$ such that the optimal rule is to accept an observation with value $a_i$ with time $\tau$ remaining if and only if $\tau\leq t_i^*$ or $i=n$. 
Set $t_0^*=0$, and $t_n^*=\infty$. For $1\leq k\leq n$ and $t\geq t_{k-1}^*$, let $V_k(t)$ denote the expected return, with time $t$ remaining, from the rule:
\begin{center}
Accept $a_i$ with time $\tau$ remaining if and only if $\tau\leq t_i^*$ or $i\geq k$.
\end{center}
Clearly, it is optimal to accept $a_k$ with time $t$ remaining if and only if $a_k\geq V_k(t)$. Thus, $t_k^*$ is the unique value of $t\geq t_{k-1}^*$ such that $V_k(t)=a_k$. For $k=1$, we have
\begin{equation*}
V_1(t)=(1-e^{-t})E_1, \qquad t\geq 0,
\end{equation*}
so that
\begin{equation*}
t_1^*=-\log(1-(a_1/E_1))=-\log(\mu_1/\mu_0)=\log(\mu_0/\mu_1).
\end{equation*}
And, inductively for $k=2,\dots,n-1$ and $t\geq t_{k-1}^*$,
\begin{align*}
V_k(t)&=\left(1-e^{-r_k(t-t_{k-1}^*)}\right)E_k+e^{-r_k(t-t_{k-1}^*)}V_{k-1}(t_{k-1}^*)\\
&=E_k-(E_k-a_{k-1})e^{-r_k(t-t_{k-1}^*)}.
\end{align*}
Thus,
\begin{equation*}
e^{-r_k(t_k^*-t_{k-1}^*)}=\frac{E_k-a_k}{E_k-a_{k-1}}=\frac{r_k(E_k-a_k)}{r_k(E_k-a_{k-1})}=\frac{\mu_k}{\mu_{k-1}},
\end{equation*}
and so
\begin{equation*}
t_k^*=t_{k-1}^*+(1/r_k)\log(\mu_{k-1}/\mu_k), \qquad k=2,\dots,n-1.
\end{equation*}
Finally, when $t_{k-1}^*\leq t\leq t_k^*$ ($k=1,2,\dots,n$),
\begin{equation}
V(t)=V_k(t)=E_k-(E_k-a_{k-1})e^{-r_k(t-t_{k-1}^*)}.
\label{eq:optimal-value}
\end{equation}

On the other hand, the prophet's value is easily computed to be
\begin{equation}
M(t)=\sum_{i=1}^n (a_i-a_{i-1})\left(1-e^{-r_i t}\right), \qquad \mbox{for all $t\geq 0$.}
\label{eq:expected-maximum}
\end{equation}

\begin{example} \label{two-point-ratio}
{\rm
Let $n=2$, and put $a_1=1$ and $a_2=K$, where $K$ is large. Let $a$ be a positive real number such that 
\begin{equation}
\log(1+t/a)<t, 
\label{eq:a-condition}
\end{equation}
and let $r_2=a/(tK)$. We will examine the ratio $R(t)=M(t)/V(t)$ as $K\to\infty$. First, by \eqref{eq:expected-maximum},
\begin{equation}
M(t)=1-e^{-t}+(K-1)\left(1-e^{-a/K}\right)\to 1-e^{-t}+a, \qquad\mbox{as $K\to\infty$}.
\label{eq:two-point-M}
\end{equation}
Next, $\mu_1=(a/tK)(K-1)\to a/t$, and $\mu_0=\mu_1+1\to (a/t)+1$, so that
\begin{equation*}
t_1^*=\log(\mu_0/\mu_1)\to\log(1+t/a), \qquad K\to\infty.
\end{equation*}
It follows that $t_1^*<t$ for sufficiently large $K$, so by \eqref{eq:optimal-value},
\begin{align}
\begin{split}
V(t)&=E_2-(E_2-a_1)e^{-r_2(t-t_1^*)}\\
&=(K-1)\left(1-e^{-a(t-t_1^*)/tK}\right)+1\\
&\to a-(a/t)\log(1+t/a)+1, \qquad K\to\infty.
\end{split}
\label{eq:two-point-V}
\end{align}
Together, \eqref{eq:two-point-M} and \eqref{eq:two-point-V} yield that
\begin{equation}
R(t)\to\frac{a+1-e^{-t}}{a+1-(a/t)\log(1+t/a)}, \qquad\mbox{as $K\to\infty$}.
\label{eq:two-point-R}
\end{equation}
In particular, if $a=1$, then \eqref{eq:a-condition} is met for every $t>0$, and \eqref{eq:two-point-R} becomes
\begin{equation}
R(t)\to\frac{2-e^{-t}}{2-\log(1+t)/t}.
\label{eq:limiting-ratio}
\end{equation}
Numerical experimentation suggests that, for the range $0<t<2$, this ratio is close to the maximum ratio over all two-valued random variables.
}
\end{example}

\begin{example} \label{three-point-ratio}
{\rm
Let $n=3$, and put $a_i=K^{i-1}$ for $i=1,2,3$, where $K$ is again assumed to be large. Let $a$ and $b$ be positive real numbers such that $a<t$ and
\begin{equation}
\log(1+a/b)<a,
\label{eq:ab-condition}
\end{equation}
and let $r_2=a/t$, and $r_3=b/(tK)$. Then
\begin{align*}
M(t)&=1-e^{-t}+(K-1)\left(1-e^{-a}\right)+K(K-1)\left(1-e^{-b/K}\right)\\
&\sim K\left(1-e^{-a}+b\right) \qquad\mbox{as $K\to\infty$.}
\end{align*}
Next, $\mu_2=(b/tK)K(K-1)=b(K-1)/t$, $\mu_1=\mu_2+(a/t)(K-1)=(a+b)(K-1)/t$, and $\mu_0=\mu_1+1$. Hence, $t_1^*=\log(\mu_0/\mu_1)=\log(1+1/\mu_1)\to 0$, and so
\begin{equation*}
t_2^*=t_1^*+(t/a)\log(1+a/b)\to (t/a)\log(1+a/b).
\end{equation*}
It follows that $t_2^*<t$ when $K$ is sufficiently large, and then
\begin{align*}
V(t)
&=K(K-1)\left(1-e^{-b(t-t_2^*)/tK}\right)+K\\
&\sim K\left[b-(b/a)\log(1+a/b)+1\right].
\end{align*}
Thus,
\begin{equation*}
R(t)\to \frac{1+b-e^{-a}}{1+b-(b/a)\log(1+a/b)}\qquad \mbox{as $K\to\infty$.}
\end{equation*}
In particular, if $a=2$ and $b=1$, then \eqref{eq:ab-condition} is satisfied, and
\begin{equation*}
R(t)\to \frac{2-e^{-2}}{2-(\log 3)/2} \doteq 1.28536, \qquad\mbox{for $t>2$}.
\end{equation*}
Note that this is the same value obtained in \eqref{eq:limiting-ratio} for $t=2$. However, by admitting three-point distributions this ratio can be achieved for {\em any} $t\geq 2$.
}
\end{example}

Observe from Table \ref{tab:alphas} that the smallest $n$ for which $a_n>1.28536$ is $n=8$. A ratio arbitrarily close to $a_8=1.29166$ can be obtained when $t\geq\log(7/\alpha_8)\doteq 3.1781$. For smaller values of $t$, however, Examples \ref{two-point-ratio} and \ref{three-point-ratio} provide larger ratios than the method discussed at the end of Section \ref{long-range}.

\begin{figure}
\begin{center}
\epsfig{file=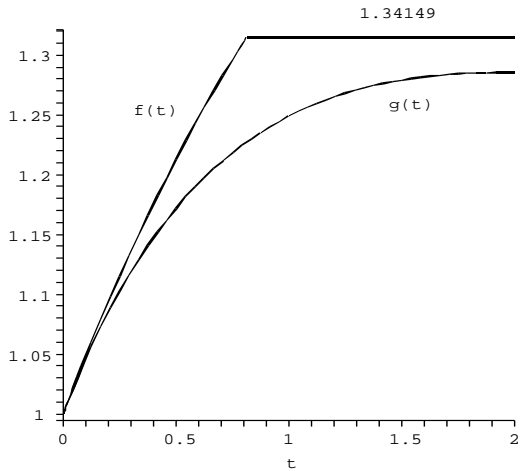,height=6.5cm,width=7.5cm}
\caption{\footnotesize The theoretical best-possible ratio bound is between the two curves $g(t)$ and $\min\{f(t),1.34149\}$, with $f(t)=2-(1-e^{-t})/t$, and $g(t)=(2-e^{-t})/\{2-\log(1+t)/t\}$.}
\label{fig:ratio-bounds}
\end{center}
\end{figure}

\bigskip
Example \ref{three-point-ratio} shows that the bound of Theorem \ref{thm:threshold-bound} is asymptotically sharp as $t\downarrow 0$ in the following sense. Let $f(t)=2-(1-e^{-t})/t$, and $g(t)=(2-e^{-t})/\{2-\log(1+t)/t\}$. That is, $g(t)$ is the right hand side of \eqref{eq:limiting-ratio}. By Example \ref{three-point-ratio}, the theoretical best-possible ratio bound is between $g(t)$ and $f(t)$. Straightforward calculations show that
\begin{equation*}
f(t)-g(t)=O(t^2)\qquad\mbox{as $t\downarrow 0$}.
\end{equation*}
This relationship is illustrated in Figure \ref{fig:ratio-bounds}, which also shows the uniform ratio bound from Theorem \ref{thm:uniform-ratio-bound}.

A similar comparison can be made for the difference bound of Corollary \ref{threshold-difference-bound}. Let $\hat{f}(t)=\beta(t)(1-e^{-t})/(\beta(t)+1-e^{-t})$. That is, $\hat{f}(t)$ is the right hand side of \eqref{eq:d-bound} for $a=0, b=1$. Let $X$ have the distribution given by $\rP(X=1)=1/(e^t+1)=1-\sP(X=(1-e^{-t})/2)$. Then $t_1^*=t$ exactly, and equations \eqref{eq:optimal-value} and \eqref{eq:expected-maximum} yield
\begin{equation}
M(t)-V(t)=\frac12\left[(1+e^{-t})(1-\exp\{-t/(e^t+1)\})-e^{-t}(1-e^{-t})\right].
\label{eq:ugly-difference}
\end{equation}
Let $\hat{g}(t)$ denote the right hand side of \eqref{eq:ugly-difference}. A straightforward calculation shows that
\begin{equation*}
\hat{f}(t)-\hat{g}(t)=O(t^3)\qquad\mbox{as $t\downarrow 0$}.
\end{equation*}
Thus, the bound of Corollary \ref{threshold-difference-bound} is asymptotically quite sharp as $t\downarrow 0$.

\section{Other renewal processes: examples in discrete time} \label{renewal}

The main purpose of this section is to show that the conclusions of Theorems \ref{thm:uniform-ratio-bound} and \ref{thm:uniform-difference-bound} may fail if the Poisson process governing the arrivals of observations is replaced by an arbitrary renewal process. The general setup is as follows. Let $T_1,T_2,\dots$ be i.i.d. random variables taking values in the positive integers, and assume that $X, X_1,X_2,\dots$ are i.i.d. nonnegative random variables, independent of the $T_i$. Call the random times $S_k=T_1+\dots+T_k$ ($k\in\NN$) the {\em renewal times}, and assume that for each $k$, the random variable $X_k$ is observed at time $S_k$. Put $S_0=0$. For $n\in\NN$, let $N_n=\max\{k:S_k\leq n\}$. In other words, $N_n$ is the number of observations that arrive by time $n$. As before, we wish to compare the values $M_n:=\sE(\max\{X_1,\dots,X_{N_n}\})$ and $V_n:=\sup_{\tau\in\TT}\sE[X_\tau \sI(S_\tau\leq n)]$, where $\TT$ is the set of all $\NN$-valued random variables (stopping rules) $\tau$ such that $\{\tau=i\}$ is measurable with respect to the $\sigma$-algebra generated by $X_1,\dots,X_i$ and $S_1,\dots,S_i$.

The problem will be easier to analyze if we represent it as follows. For each $j\in\NN$, define
\begin{equation*}
Y_j=\begin{cases}
X_k, & \mbox{if $j=S_k \quad (k=1,2,\dots)$},\\
0, & \mbox{otherwise}.
\end{cases}
\end{equation*}
It is not difficult to see that $M_n=\sE(Y_1\vee\dots\vee Y_n)$, and $V_n=\sup\{\sE Y_\tau: \tau$ is a stopping rule for $Y_1,\dots,Y_n\}$. Thus, the problem is reduced to that of stopping an ordinary sequence of random variables, and standard methods can be applied to solve it.

Observe that the $Y_j$ are, in general, neither independent nor identically distributed.  However, there is one important exception.

\begin{example} \label{binomial-process}
{\rm
Let $0<p<1$, and assume that $\rP(T_1=k)=(1-p)^{k-1}p$ for $k=1,2,\dots$. This yields the {\em binomial process}, which has the property that the events $\{j$ {\em is a renewal time}$\}$, $j\in\NN$, are mutually independent and have probability $p$. Since the $X_i$ are i.i.d., this implies that the $Y_j$ are i.i.d. with common distribution $pF+(1-p)\delta_{\{0\}}$, where $F$ is the distribution of $X$, and $\delta_{\{0\}}$ denotes Dirac measure at zero. It follows from Theorem~A and the above representations that $M_n\leq a_n V_n$. Similarly, if $X$ is $[0,1]$-valued then so is $Y_1$, and Theorem~B implies that $M_n-V_n\leq b_n$. The sharpness of these inequalities depends on the value of $p$: the first bound is sharp if $1-p\leq (\eta_{0,n}(\alpha_n))^{1/n}$, that is, if $p\geq (\alpha_n/(n-1))^{1/n}$. Likewise, the second bound is attained if $p\geq (\beta_n/(n-1))^{1/n}$. It is not clear how sharp the bounds are when $p$ is smaller than the indicated values.
}
\end{example}

The next example shows that the best possible ratio and difference bounds are, in general, strictly greater than $a_n$ and $b_n$.

\begin{example} \label{counter-example}
{\rm
Fix $n\in\NN$. Let $0<p<1$, and assume that $\rP(T_1=1)=p=1-\sP(T_1=n)$. We compute $V_n$ by backward induction. For $i=1,2,\dots,n$, let $\gamma_i$ denote the supremum, over all stopping times $\tau$ such that $i\leq\tau\leq n$, of $\rE[Y_\tau|\,i$ {\em is a renewal time}$]$. Then $\gamma_n=\sE X$, and $\gamma_i=\sE(X\vee p\gamma_{i+1})$ for $i=1,2,\dots,n-1$, since if $i$ is a renewal moment, then the next renewal moment is either $i+1$ or $i+n$, and $i+n$ is beyond the time horizon. Finally,
\begin{equation}
V_n=p\gamma_1+(1-p)\gamma_n.
\label{eq:Vn}
\end{equation}
Now let $X$ have a distribution on two points $\eps$ and $1$, where $0<\eps<1$, and the probability $\pi:=\sP(X=1)$ is chosen so that
\begin{equation}
p\sE X=\eps.
\label{eq:mean-condition}
\end{equation}
It follows immediately that $\gamma_i=\sE X$ for $i=1,\dots,n$, and hence, by \eqref{eq:Vn}, $V_n=\sE X$. On the other hand,
\begin{align*}
M_n&=\sE X_1+\sum_{k=2}^n \sE(X_2\vee\dots\vee X_k-X_1)^+\sP(N_n=k)\\
&=\sE X+\sum_{k=2}^{n-1}p^k(1-p)(1-\eps)(1-\pi)\{1-(1-\pi)^{k-1}\}\\
&\qquad\qquad\qquad\qquad\qquad +p^n(1-\eps)(1-\pi)\{1-(1-\pi)^n\}\\
&=\sE X+(1-\eps)p\pi \left[\frac{1-\{p(1-\pi)\}^n}{1-p(1-\pi)}-1\right],
\end{align*}
where the last equality follows after routine simplification. Since $\rE X=\pi+(1-\pi)\eps$, \eqref{eq:mean-condition} implies that $\eps=p\pi/\{1-p(1-\pi)\}$, and $1-\eps=(1-p)/\{1-p(1-\pi)\}$. Thus, we obtain the expressions
\begin{equation*}
D_n:=M_n-V_n=\frac{p(1-p)\pi}{1-p(1-\pi)}\left[\frac{1-\{p(1-\pi)\}^n}{1-p(1-\pi)}-1\right]
\end{equation*}
and
\begin{equation*}
R_n:=\frac{M_n}{V_n}=1+p(1-p)\left[\frac{1-\{p(1-\pi)\}^n}{1-p(1-\pi)}-1\right].
\end{equation*}
Now as $\pi\downarrow 0$, $R_n$ increases to $1+p^2-p^{n+1}$, which is maximized for $p=(2/(n+1))^{1/(n-1)}$. It follows that $R_n$ can be arbitrarily close to
\begin{equation*}
c_n:=1+\left(\frac{2}{n+1}\right)^{2/(n-1)}\left(\frac{n-1}{n+1}\right).
\end{equation*}
Observe that $\lim_{n\to\infty}c_n=2$. Thus, the conclusion of Theorem \ref{thm:uniform-ratio-bound} fails to hold for this case when $n$ is sufficiently large. (In fact, $c_5\doteq 1.3849>1+\alpha_0$.)

As for the difference $D_n$, note that
\begin{equation}
\lim_{n\to\infty}D_n=\frac{p^2(1-p)\pi(1-\pi)}{\{1-p(1-\pi)\}^2}.
\label{eq:Dn}
\end{equation}
For fixed $p$, this is maximized at $\pi=(1-p)/(2-p)$, and substituting this into \eqref{eq:Dn} yields that $\lim_{n\to\infty}D_n=p^2/4$. Thus, if we choose $n$ sufficiently large, $p$ sufficiently close to $1$, and $\pi=(1-p)/(2-p)$, then $D_n$ will be arbitrarily close to $1/4$.

Note that it is not known whether $\limsup_{n\to\infty}b_n<1/4$, though Table \ref{tab:alphas} suggests this should be the case. If this is true, then the conclusion of Theorem \ref{thm:uniform-difference-bound} too fails to hold for this example.
}
\end{example}

The last example raises an interesting question: do there exist a renewal process (in discrete or continuous time) and a distribution for $X$ such that $M(t;X)/V(t;X)>2$, or (if $X$ is $[0,1]$-valued) $M(t;X)-V(t;X)>1/4$? If not, why do these classical constants for the independent case appear as upper bounds in a problem concerning i.i.d. random variables? These questions may be addressed in a future paper.

\end{document}